\documentclass[12pt]{article}
\usepackage{amssymb}

\newcommand{\SUM}{\raisebox{-0.4ex}{\mbox{\Large $\Sigma$}}}

\newcommand{\EE}{{\mathbb E}}
\newcommand{\F}{{\mathbb F}}
\newcommand{\LL}{\Lambda_3}
\newcommand{\PP}{\mathbb{P}}

\newtheorem{theorem}{Theorem}
\newtheorem{lemma}{Lemma}

\newtheorem{proposition}{Proposition}

\title{On the Structure of Sets with Few Three-Term Arithmetic Progressions}
\author{Ernie Croot}

\begin{document}

\maketitle

\section{Introduction}

Given a function $f : \F_{p^n} \to [0,1]$, and a subset $W \subseteq \F_{p^n}$, 
we define
$$
\EE(f | W)\ =\ |W|^{-1} \SUM_{m \in W} f(m).
$$
If no set $W$ is given, then we just assume $W = \F_{p^n}$, and then we get
$$
\EE(f)\ =\ \EE(f | \F_{p^n})\ =\ p^{-n} \SUM_{m \in \F_{p^n}} f(m).
$$
Define
$$
\LL(f)\ =\ p^{-2n} \SUM_{m,d} f(m)f(m+d)f(m+2d).
$$
In the case where $f$ is an indicator function for some set $S \subseteq \F_{p^n}$,
we have that $\LL(f)$ is the normalzed count of the number of three-term 
arithmetic progressions $m,m+d,m+2d \in S$.  Note that $\LL(f) \geq 0$, unless
$\EE(f) = 0$, because of the contribution of trivial progressions where
$d = 0$.
\bigskip

Of central importance to the subject of additive combinatorics is that of
determining when a subset of the integers $\{1,...,N\}$ contains a 
$k$-term arithmetic progression.  This subject has a long history, and
we will not mention it here; however, the specific problem in this area
which motivated our paper, and which is due to B. Green \cite{AIM}, 
is as follows:
\bigskip

\noindent {\bf Problem.}  Given $0 < \alpha \leq 1$, suppose 
$S \subseteq {\mathbb F}_p$ satisfies $|S| \geq \alpha p$, and
has the least number of three-term arithmetic progressions.  
What is $\LL(S)$ ?
\bigskip

It seems that the only hope of answering a question like this is to 
understand the structure of these sets $S$.  In this paper we address 
the analogous problem in $\F_{p^n}$, where $p$ and $\alpha$ are held
fixed, while $n$ tends to infinity.  The results we prove
are not of a type that would allow us to dedcue $\LL(S)$, but they do
reveal that these sets $S$ are very highly structured.  Such results
can perhaps be deduced from the work of B. Green \cite{green}, which
makes use of the Szemer\'edi regularity lemma, but our theorems below
are proved using basic harmonic analysis.

\begin{theorem} \label{main_theorem1}  Let $0 < \alpha \leq 1$.  
Suppose that $S$ is a subset of 
$\F_{p^n}$, such that $\LL(S)$ is minimal, subject to the constraint
$$
|S|\ \geq\ \alpha p^n.
$$
Then, there exists a subgroup (or subspace) 
$$
W\ \leq\ \F_{p^n},\ {\rm dim}(W)\ =\ n-o(n), 
$$
such that $S$ is approximately a union of $p^{o(n)}$ cosets 
of $W$; more precisely, there is a set $A$ of size $p^{o(n)}$ such that
$$
|S\ \Delta\ A+W|\ =\ o(p^n). \footnote{The notation $B \Delta C$ means
the symmetric difference between $B$ and $C$.}
$$
\end{theorem}

Our second theorem is a slighly more abstract version of this one, where
instead of sets $S$, we have a function $f : \F_{p^n} \to [0,1]$.  

\begin{theorem} \label{main_theorem2}  Let $0 < \alpha \leq 1$.  
Suppose that 
$$
f\ :\  \F_{p^n} \to [0,1]
$$
such that $\LL(f)$ is minimal, 
subject to the constraint that 
$$
\EE(f)\ \geq\ \alpha\ >\ 0.
$$  
Then, there exists a subgroup $W \subseteq \F_{p^n}$ of 
dimension $n - o(n)$, such that 
$f$ is approximately an indicator function on cosets of $W$, in the
following sense:  There is a function 
$$
h\ :\ \F_{p^n}\ \to\ \{0,1\},
$$
which is constant on cosets of $W$ (which means $h(a) = h(a+w)$ for
all $w \in W$), such that 
$$
\EE(|f(m) - h(m)|)\ =\ o(1).
$$
\end{theorem}

It would seem that Theorem \ref{main_theorem1} is a corollary of
Theorem \ref{main_theorem2}; however, with a little thought one sees
this is not the case.  Nonetheless, we will prove a third theorem,
from which we will deduce both 
Theorem \ref{main_theorem1} and Theorem \ref{main_theorem2}.

\section{Proofs}

\subsection{Additional Notation}

We will require a little more notation.

Given any three subsets $U,V,W \subseteq \F_{p^n}$, define  
$$
T_3(f|U,V,W)\ =\ \SUM_{m\in U, m+d \in V, m+2d \in W} f(m)f(m+d)f(m+2d).
$$
We note that this implies $T_3(1|U,U,U)$ is the number of three-term 
progressions belonging to a set $U$.
\bigskip

Given a subspace $W$ of $\F_{p^n}$, and given a function 
$$
f\ :\ \F_{p^n}\ \to\ [0,1],
$$
we define
$$
f_W(m)\ =\ {1 \over |W|} (f*W)(m)\ =\ {1 \over |W|} \SUM_{w \in W} f(m+w).
$$
This function has a number of properties:  First, we note that 
$f_W(m)$ is constant on cosets of $W$, in the sense that 
$$
{\rm for\ all\ }w \in W,\ f_W(m)\ =\ f_W(m+w).
$$
Thus, it makes sense to write 
$$
f_W(m + W)\ =\ f_W(m).
$$
We also have that 
\begin{equation} \label{fw_expectation}
\EE(f_W)\ =\ \EE(f).
\end{equation}
Finally, if $V$ is the orthogonal complement of $W$ (with respect to
the standard basis), then 
\begin{equation} \label{fw_fourier}
{\rm if\ }v \in V,\ {\rm then\ }\hat f_W(v) = \hat f(a);\ {\rm and,}\ 
{\rm if\ }v \not \in V,\ {\rm then\ }\hat f_W(v) = 0.
\end{equation}
\bigskip

We will also define the $L^2$ norm of a function $f : \F_{p^n} \to 
{\mathbb C}$ to be 
$$
||f||_2\ =\ \Bigl ( p^{-n} \SUM_m |f(m)|^2 \Bigr )^{1/2}.
$$

\subsection{Theorem \ref{main_theorem3}, and Proofs of
Theorems \ref{main_theorem1} and \ref{main_theorem2}}

Theorems \ref{main_theorem1} and \ref{main_theorem2} are corollaries 
of the following theorem:

\begin{theorem} \label{main_theorem3}
Let $\epsilon > 0$, and suppose that 
$$
f\ :\ \F_{p^n}\ \to\ [0,1]
$$
has the following property:  For every subspace $W$ of $\F_{p^n}$ of
codimension at most $\Delta^{-2}$, where
$$
\Delta\ =\ (\epsilon^6/2^{13}p^2) \exp(-16\epsilon^{-1}c_p \log p),
$$
where $c_p$ is a certain constant appearing in Theorem \ref{meshulam_theorem}
below, suppose that 
$$
\EE(|f(m)-f_W(m)|)\ >\ \epsilon.
$$

Then, there exists a function 
$$
g\ :\ \F_{p^n}\ \to\ [0,1]
$$
such that 
$$
\EE(g)\ =\ \EE(f),\ {\rm and\ } \LL(g)\ <\ \LL(f) - \Delta.
$$
\end{theorem}

\noindent {\bf Comment.}  Using the Lemma \ref{random_lemma} below 
we can deduce the stronger conclusion that there exists 
$$
g\ :\ \F_{p^n}\ \to\ \{0,1\}
$$
(so, $g$ is an indicator function) such that 
\begin{equation} \label{g_good}
\EE(g)\ \geq\ \EE(f),\ {\rm and\ } \LL(g)\ <\ \LL(f) - \Delta + O(p^{-n/3}).
\end{equation}
\bigskip

\begin{lemma}  \label{random_lemma}
Suppose that $j\ :\ \F_{p^n}\ \to\ [0,1]$.  There exists an 
indicator function $j_2\ :\ \F_{p^n}\ \to\ \{0,1\}$, such that 
$$\\
\EE(j_2)\ \geq\ \EE(j),\ \LL(j_2)\ =\ \LL(j) + O(p^{-n/3}),
$$\\
and such that for every subspace $W$ of codimension at most 
$n^{1/2}$ we have\footnote{The codimension $n^{1/2}$ condition
can be improved; however, it is good enough for our purposes, and
it is larger than $\Delta^{-2}$, where $\epsilon = 1/\log\log n$,
as will appear in later applications.} 
that for every $m \in \F_{p^n}$,
$$
(j_2)_W(m)\ =\ j_W(m) + O(1/n).
$$
\end{lemma}

In order to prove this lemma we will need to use a theorem of
Hoeffding (see \cite{hoeffding} or \cite[Theorem 5.7]{mcdiarmid})

\begin{proposition} \label{hoeffding_prop} 
Suppose that $z_1,...,z_r$ are independent real random variables 
with $|z_i| \leq 1$.  Let $\mu = \EE(z_1 + \cdots + z_r)$,
and let $\Sigma = z_1 + \cdots + z_r$.  Then,
$$\\
\PP(|\Sigma - \mu| > rt)\ \leq\ 2 \exp(-rt^2/2).
$$
\end{proposition}

\noindent {\bf Proof of the Lemma.}
The proof of this lemma is standard:  Given $j$ as in the theorem above,
let $j_0$ be a random function from $\F_{p^n}$ to $\{0,1\}$, where 
$j_0(m) = 1$ with probability $j(m)$, and equals $0$ with probability
$1-j(m)$; moreover, $j_0(m)$ is indepedent of all the other
$j_0(m')$.  Then, one can easily show that with probability $1-o(1)$,
\begin{equation} \label{j11}
p^{-n} \SUM_m j_0(m)\ =\ \EE(j) + O(p^{-n/3}),\ {\rm and\ }
\LL(j_0)\ =\ \LL(j) + O(p^{-n/3}).
\end{equation}
Furthermore, we claim that with probability $1-o(1)$ we will have that for any 
subspace $W$ of codimension at most $n^{1/2}$, 
\begin{equation} \label{W_hope}
(j_0)_W(m)\ =\ j_W(m) + O(1/n).
\end{equation}
This can be seen as follows:   For
a fixed $W$ we need an upper bound on the probability that 
$$
|(j_0)_W(m) - j_W(m)|\ >\ 1/n.
$$
This is the same as showing
$$
|\Sigma|\ >\ |W|/n,
$$
where 
$$
\Sigma\ =\ \SUM_{w \in W} z_w(m),\ {\rm where\ } 
z_w(m)\ =\ j_0(m+w) - j(m+w).
$$
Note that all the $z_w$ are independent
and satisfy $|z_w| \leq 1$ and $\EE(z_w) = 0$.  So, from 
Proposition \ref{hoeffding_prop} we deduce that 
$$
\PP(|\Sigma| > |W|/n)\ \leq\ 2\exp(-|W|/2n^2).
$$
Now, since the number of such subspaces $W$ is at 
most the number of sequences of $n^{1/2}$ possible basis vectors, 
which is $O(p^{n^{3/2}})$, we deduce that the probability that there exists 
a subspace $W$ of codimension at most $n^{1/2}$ satisfying
$$
|(j_0)_W(m) - j_W(m)|\ >\ 1/n
$$\\
is $O(p^{n^{3/2}} \exp(-|W|/2n^2)) = o(1)$.  Thus, (\ref{W_hope}) holds for
all such $W$ with probability $1-o(1)$ (in fact, the explicit constant in the
$O(1)$ can be taken to be $1$ once $n$ is sufficiently large).  
\bigskip

We deduce now that there is an instantiation of $j_0$, call it $j_1$, such that 
both (\ref{j11}) and (\ref{W_hope}) hold.  Then, by reassigning at most 
$O(p^{2n/3})$ places $m$ where $j_1(m) = 0$ to the value $1$, or 
from the value $0$ to the value $1$, we arrive at a function $j_2$
having the claimed propertes of the lemma. \hfill $\blacksquare$
\bigskip

\noindent {\bf Proof of Theorem \ref{main_theorem1}.}
To prove Theorem \ref{main_theorem1}, we begin by letting $f$ be 
the indicator function for the set $S$, and we let 
$$
\epsilon\ =\ {1 \over \log\log n}.
$$

Now suppose that 
\begin{equation} \label{assumption1}
\EE(|f(m) - f_W(m)|)\ \leq\ \epsilon,
\end{equation}
for some subspace $W$ of codimension at most $\Delta^{-2}$.  
Let $h(m)$ be $f_W(m)$ rounded to the nearest integer.  Clearly,
$h(m)$ is constant on cosets of $W$, and from the fact that
$$
|h(m) - f_W(m)|\ \leq\ |f(m) - f_W(m)|,
$$
we deduce that 
\begin{eqnarray}
\EE(|f(m) - h(m)|)\ &\leq&\ \EE(|h(m) - f_W(m)|) + \EE(|f(m) - f_W(m)|)\nonumber \\
&\leq&\ 2 \EE(|f(m) - f_W(m)|)\nonumber \\
&\leq&\ 2\epsilon. \nonumber
\end{eqnarray}
But since $h$ is constant on cosets of $W$, and only assumes the values
$0$ or $1$, we deduce that $h$ is the indicator function for some 
set of the form $A+W$.  Thus, we deduce
$$
|S\ \Delta\ A+W|\ \leq\ 2 \epsilon p^n,
$$
where $W$ has dimension $n-o(n)$.  This then proves Theorem \ref{main_theorem1}
under the assumption (\ref{assumption1}).  

Next, suppose that 
\begin{equation} \label{assumption2}
\EE(|f(m) - f_W(m)|)\ >\ \epsilon.
\end{equation}
for every subspace $W$ of codimension at most $\Delta^{-2}$.
Then, from the comment following 
Theorem \ref{main_theorem3}, there exists an indicator 
function $g$ satisfying (\ref{g_good}).  If we let
$S'$ be the set for which $g$ is an indicator function, then
one sees that $S'$ has fewer three-term arithmetic progressions
than does $S$, while $\EE(S') \geq \EE(S)$.  This is a contradiction,
and thus the theorem is proved.
\hfill $\blacksquare$
\bigskip  

\noindent {\bf Proof of Theorem \ref{main_theorem2}.}
Let $j(m) = f(m)$, and then let 
$$
\ell(m)\ =\ j_2(m)\ :\ \F_{p^n}\ \to\ \{0,1\},
$$ 
where $j_2(m)$ is as given in Lemma \ref{random_lemma}.  
Note that this implies that
$$
\EE(\ell)\ \geq\ \EE(f),\ \LL(\ell)\ =\ \LL(f) + O(p^{-n/3}),
$$
and that for any subspace $W$ of codimension at most $n^{1/2}$,
\begin{equation} \label{ell_condition2}
\ell_W(m)\ =\ f_W(m) + O(1/n).
\end{equation}
\bigskip

Next let 
$$
\epsilon\ =\ {1 \over \log\log n},
$$
and suppose that there exists a subspace $W$ of codimension at 
most $\Delta^{-2}$ such that 
\begin{equation} \label{subspace_assume}
\EE(|\ell(m) - \ell_W(m)|)\ \leq\ \epsilon.
\end{equation}
Then, if we let $h(m)$ equal $f_W(m)$ rounded to the nearest integer,
we will have from (\ref{ell_condition2}) that
\begin{eqnarray} \label{jfW}
\EE(|h(m) - f_W(m)|)\ &\leq&\ \EE(|\ell(m) - f_W(m)|)\nonumber \\
&\leq&\ \EE(|\ell(m) - \ell_W(m)|) + O(1/n) \nonumber \\
&\leq&\ \epsilon + O(1/n). 
\end{eqnarray}

Let $V$ be the orthogonal complement of $W$.  From (\ref{jfW}) we
know that at most
$$
(\epsilon^{1/2} + O(\epsilon^{-1/2}/n)) |V|
$$ 
values $v \in V$ satisfy
$$
|h(v) - f_W(v)|\ \geq\ \epsilon^{1/2}.
$$
Let $V' \subseteq V$ be those $v \in V$ satisfying the
reverse inequality
$$
|h(v) - f_W(v)|\ <\ \epsilon^{1/2}.
$$

Suppose $v \in V'$ and $h(v) = 0$.  Then, 
$f_W(v) < \epsilon^{1/2}$, and we have 
\begin{equation} \label{j0}
\SUM_{m \in v+W} |f(m) - h(m)|\ =\ |W|f_W(v)\ <\ 
|W|\epsilon^{1/2}.
\end{equation}
On the other hand, if $v \in V'$ and $h(v) = 1$, then 
$f_W(v) > 1 - \epsilon^{1/2}$,
and so
\begin{equation} \label{j1}
\SUM_{m \in v+W} |f(m) - h(m)|\ =\ |W|(1 - f_W(v))\ <\ |W|\epsilon^{1/2}.
\end{equation}

Combining (\ref{j0}) with (\ref{j1}) we deduce that 
\begin{eqnarray}
\EE(|f(m) - h(m)|)\ &\leq&\ \epsilon^{1/2} + (|V| - |V'|)|V|^{-1}
\nonumber \\ 
&\leq&\ 2\epsilon^{1/2} + O(\epsilon^{-1/2}/n).
\end{eqnarray}
Our theorem is now proved in this case (assuming there exists
a subspace $W$ satisfying (\ref{subspace_assume}) ).
\bigskip

To complete the proof, we will assume that there are no
subspaces of codimension at most $\Delta^{-2}$ satisfying
(\ref{subspace_assume}).  Since $\ell$ then satisfies
the hypotheses of Theorem \ref{main_theorem3}, we deduce from
Theorem \ref{main_theorem3} that there exists a function
$g : \F_{p^n} \to [0,1]$ such that 
$$
\EE(g)\ =\ \EE(\ell)\ \geq\ \EE(f) \geq \alpha,
$$
and
$$
\LL(g)\ <\ \LL(\ell) - \Delta\ =\ \LL(f) - \Delta + O(p^{-n/3}).
$$
This then contradicts the fact that $\LL(f)$ was minimal,
given $\EE(f) \geq \alpha$.  Our theorem is now proved.
\hfill $\blacksquare$

\section{Proof of Theorem \ref{main_theorem3}}

Let $\Delta$ be as in the statement of Theorem \ref{main_theorem3}.

As is well-known,
$$
\LL(f)\ =\ p^{-3n} \SUM_{a \in \F_{p^n}} \hat f(a)^2 \hat f(-2a).
$$
If we let $A$ denote the set of all $a \in \F_{p^n}$ where 
$$
|\hat f(a)|\ >\ \Delta p^n,
$$
then we clearly have
\begin{equation} \label{LLf}
\LL(f)\ =\ p^{-3n} \SUM_{a \in A} \hat f(a)^2 \hat f(-2a)\ +\ E,
\end{equation}
where
\begin{equation} \label{LLfE}
|E|\ \leq\ \Delta p^{-n} ||\hat f||_2^2\ \leq\ \Delta.
\end{equation}
A simple application of Parseval's identity also shows that $|A|$ is small:  
We have
$$
|A| \Delta^2 p^{2n}\ \leq\ p^n ||\hat f||_2^2\ \leq\ p^{2n},
$$
which implies
$$
|A|\ \leq\ \Delta^{-2}. 
$$
\bigskip

Let $V$ be the additive subgroup of $\F_{p^n}$ generated by the elements of $A$,
and let $W$ be the orthogonal complement of $V$; that is,

$$
W\ =\ \{w \in \F_{p^n}\ :\ {\rm for\ every\ } v \in V,\ w\cdot v\ =\ 0\}.
\footnote{The product $w\cdot v$ here denotes the dot product with respect 
to the standard basis of the vector space $\F_{p^n}$, not the product defined 
for the multiplicative structure of $\F_{p^n}$. } 
$$
From (\ref{LLf}), (\ref{LLfE}), and (\ref{fw_fourier}) we deduce that 
\begin{equation} \label{fWf}
\LL(f_W)\ \leq\ \LL(f) + \Delta.
\end{equation}

Since $W$ is an additive subgroup of $\F_{p^n}$, we will use the 
standard representation for the cosets of $W$, given by
$$
v + W,\ {\rm where\ } v \in V.
$$
This canonical representation for
the cosets of $W$ has the following important property.

\begin{lemma} \label{orthogonal_lemma}  Suppose that 
$h : \F_{p^n} \to [0,1]$.  Then,
$$
T_3(h)\ =\ \SUM_{v_1,v_2,v_3 \in V \atop v_1 + v_3 = 2v_2} 
T_3(h|v_1 + W, v_2 + W, v_3 + W).
$$ 
\end{lemma}

\noindent {\bf Proof.}  The lemma will follow if we can just show that 
$v_1 + w_1, v_2 + w_2, v_3 + w_3$, $v_1,v_2,v_3 \in V$ and $w_1,w_2,w_3 \in W$,
are in arithmetic progression implies $v_1,v_2,v_3$ are in arithmetic progression:
If
$$
(v_1 + w_1) + (v_3 + w_3)\ =\ 2(v_2 + w_2),
$$
then
$$
v_1 + v_3 - 2v_2\ =\ -w_1 - w_3 + 2w_2.
$$
Now, as $V \cap W = \{0\}$, we deduce that
$$
v_1 + v_3 - 2v_2\ =\ 0,
$$
whence $v_1,v_2,v_3$ are in arithmetic progression. \hfill $\blacksquare$
\bigskip

Now let 
\begin{equation} \label{sflat}
V'\ :=\ \{ v \in V\ :\ f_W(v + W)\ \in\ [\epsilon/4,\ 1- \epsilon/4]\};
\end{equation}
that is, these cosets are all the places where $f_W$ is not 
``too close'' to being an indicator function.

\subsection{Construction of the Function $g$}

To construct the function $g$ with the properties claimed by our Theorem,
we start with the following lemma:

\begin{lemma}  Suppose $h_1 : \F_{p^n} \to [0,1]$, let $\beta = \EE(h_1)$, and  
let $h_2(n) = 1 - h_1(n)$.  Then,
$$
\LL(h_1) + \LL(h_2)\ =\ 1 - 3\beta + 3\beta^2.
$$
\end{lemma}
\bigskip

\noindent {\bf Proof.}  We first realize that for $a \neq 0$, 
$\hat h_1(a) = -\hat h_2(a)$.  Thus,
\begin{eqnarray}
\LL(h_1) + \LL(h_2)\ &=&\ 
p^{-3n} \SUM_a (\hat h_1(a)^2 \hat h_1(-2a) + \hat h_2(a)^2 \hat h_2(-2a))
\nonumber \\
&=&\ p^{-3n}(\hat h_1(0)^3 + \hat h_2(0)^3) \nonumber \\ 
&=&\ \beta^3 + (1-\beta)^3.\ \ \ \ \ \ \ \ \blacksquare \nonumber
\end{eqnarray}
\bigskip

Now, let $\ell$ be the unique integer satisfying
$$
4/\epsilon\ \leq\ p^\ell\ <\ 4p/\epsilon,
$$
and let $S$ be any subspace of $W$ of codimension
$\ell$.  Let $T$ be the complement of $S$ relative to $W$ (not 
{\it orthogonal} complement, as we have used earlier), and set 
$$
\beta\ =\ {|T| \over |W|}\ =\ {|W| - |S| \over |W|}\ =\ 
1 - p^{-\ell}\ \geq\ 1 - \epsilon/4,
$$
which is the density of $T$ relative to $W$.  
Then, from the above lemma, we deduce that
$$
T_3(S) + T_3(T)\ =\ (1 - 3\beta + 3\beta^2)|W|^2,
$$
$T_3(S)$ clearly equals $(1-\beta)^2 |W|^2$, because given any pair
of elements $m,m+d \in S$, since $S$ is a subspace we also must have 
$m+2d \in S$; and, note that there are $(1-\beta)^2 |W|^2$ ordered pairs 
$m,m+d$ in $S$.  Thus, we deduce
$$
T_3(T)\ =\ (2\beta^2 - \beta)|W|^2.
$$
We also have that if $b_1 + W, b_2 + W, b_3 + W$ are cosets that are in arithmetic
progression, in the sense that there is a triple $m,m+d,m+2d$, belonging to 
$b_1 + W, b_2 + W,$ and $b_3 + W$, respectively, then 
$$
T_3(1| b_1 + T, b_2 + T, b_3 + T)\ =\ (2\beta^2 - \beta)|W|^2.
$$
\bigskip

We now define the function $g : \F_{p^n} \to [0,1]$ as follows:  
Given $v \in V, w \in W$,
we have
$$
g(v+w)\ =\ \left \{\begin{array}{rl} f_W(v),\ &{\rm if\ } 
v \not \in V'; \\
\beta^{-1} T(w)f_W(v),\ &{\rm if\ } v \in V'. 
\end{array}\right.
$$

It is easy to see that 
$$
\EE(g)\ =\ \EE(f_W)\ =\ \EE(f);
$$
We also observe, from Lemma \ref{orthogonal_lemma}, that 
$$
T_3(g)\ =\ \SUM_{v_1,v_2,v_3 \in V \atop v_1 + v_3=2v_2} 
T_3(g | v_1 + W, v_2 + W, v_3 + W).
$$
This sum has eight types of terms, according to whether each of 
$v_1,v_2,v_3$ lie in $V'$ or not.  

First, consider the case where all of 
\begin{equation} \label{main_case}
v_1, v_2,v_3\ \in\ V'.
\end{equation}
In this case we have
\begin{eqnarray}
T_3(g | v_1 + W, v_2 + W, v_3 + W)\ &=&\ 
\beta^{-3} f_W(v_1)f_W(v_2)f_W(v_3) T_3(T) \nonumber \\
&=&\ f_W(v_1)f_W(v_2)f_W(v_3)|W|^2(2\beta^{-1} - \beta^{-2}) \nonumber \\
&\leq&\ f_W(v_1)f_W(v_2)f_W(v_3) |W|^2(1 - p^{-2\ell}) \nonumber \\
&<&\ f_W(v_1)f_W(v_2)f_W(v_3) |W|^2(1 - \epsilon^2/16p^2). \nonumber
\end{eqnarray}
This last inequality follows from the fact that 
$$
p^\ell\ <\ 4p / \epsilon.  
$$

Now, as
$$
T_3(f_W | v_1 + W, v_2 + W, v_3 + W)\ =\ f_W(v_1)f_W(v_2)f_W(v_3) |W|^2,
$$
we deduce that if (\ref{main_case}) holds, then
$$
T_3(g | v_1 + W, v_2 + W, v_3 + W)\ \leq\ T_3(f_W | v_1 + W, v_2 + W, v_3 + W)
(1-\epsilon^2/16p^2).
$$
\bigskip

On the other hand, if any of $v_1,v_2,v_3$ fail to lie in $V'$, 
then we will get that 
$$
T_3(g | v_1 + W, v_2 + W, v_3 + W)\ =\ T_3(f_W | v_1 + W, v_2 + W, v_3 + W).
$$
To see this, consider all the cases where 
$v_1$ fails to lie in $V'$.  In this case, we clearly have 
\begin{eqnarray}
T_3(g | v_1 + W,v_2 + W, v_3 + W)\ &=&\ 
\SUM_{m_1 \in v_2 + W, m_2 \in v_3 + W} 
f_W(v_1)g(m_1)g(m_2)\nonumber \\
&=&\ f_W(v_1) (|W|^2 f_W(v_2) f_W(v_3)) \nonumber \\
&=&\ T_3(f_W | v_1 + W, v_2 + W, v_3 + W). \nonumber
\end{eqnarray}
The cases where $v_2$ or $v_3$ fail to lie in 
$V'$ are identical to this one.

Putting together the above observations we deduce that
\begin{eqnarray} \label{T3f2}
T_3(g)\ &\leq&\ T_3(f_W) - 
(\epsilon^2/16p^2) \SUM_{v_1,v_2,v_3 \in V' \atop v_1 + v_3 = 2v_2}  
T_3(f_W | v_1 + W, v_2 + W,v_3 + W) \nonumber \\
&\leq&\ T_3(f_W) - (\epsilon^5/1024p^2)|W|^2 T_3(V'). 
\end{eqnarray}
This last inequality follows from the fact that $f_W(v) \geq \epsilon/4$
for $v \in V'$.

\subsection{A Lower Bound for $|V'|$}

In order to give a lower bound for $T_3(V')$, we will first
need a lower bound for $|V'|$.

We begin by noting that if $v$ belongs
to $V$, but not $V'$, then either $f_W(v) < \epsilon/4$ or 
$f_W(v) > 1 - \epsilon/4$.  Suppose the former holds.  Then,
we have 
\begin{eqnarray} \label{v`1}
\SUM_{m \in v+W} |f(m) - f_W(m)|\ \leq\ |W|f_W(v) + \sum_{m \in v+W} 
f(m)\ &=&\ 2|W|f_W(v)\nonumber \\
&<&\ \epsilon |W|/2. 
\end{eqnarray}
On the other hand, if $f_W(v) > 1-\epsilon/4$, then we have
\begin{eqnarray} \label{v`2}
\SUM_{m \in v+W} |f(m) - f_W(m)|\ &\leq&\ \SUM_{m \in v+V} (1 - f(m))
+ \SUM_{m \in v+W} (1 - f_W(m))\nonumber \\ 
&=&\ 2|W| - 2|W|f_W(v) \nonumber \\
&<&\ \epsilon|W|/2. 
\end{eqnarray}
Putting together (\ref{v`1}) and (\ref{v`2}) we deduce that
$$
\SUM_{v\in V\setminus V'} \SUM_{m \in v+W} |f(m) - f_W(m)|\ < \epsilon |W|(|V|-|V'|)/2.
$$
We also have the trivial upper bound
$$
\SUM_{v \in V'} \SUM_{m \in v+M} |f(m) - f_W(m)|\ \leq\ 
|W| |V'|.
$$
Thus, 
$$
|V|^{-1} (|V'| + \epsilon(|V| - |V'|)/2)\ >\ \EE(|f(m) - f_W(m)|)\ >\ 
\epsilon.
$$ 
(The second inequality is one of the hypotheses of the Theorem.)
It follows that
\begin{equation} \label{V'_lower}
|V'|\ >\ {\epsilon |V| \over 2(1 - \epsilon/2)}\ >\ \epsilon |V|/2.
\end{equation}

\subsection{Some Results of Meshulam and Varnavides}

Using our lower bound for $|V'|$, we will need the following
result of Meshulam \cite{meshulam} to obtain a lower bound for $T_3(V')$:

\begin{theorem} \label{meshulam_theorem}  
Suppose that $S \subseteq \F_{p^n}$ satisfies 
$|S| \geq c_p p^n/n$, where $c_p > 0$ is a certain constant depending only on $p$.
Then, $S$ contains a non-trivial three-term arithmetic progression.
\end{theorem}

If we combine this with an idea of Varnavides \cite{varnavides}, we 
get the following theorem.

\begin{theorem} \label{meshulam2} 
Suppose that $S \subseteq \F_{p^n}$ satisfies 
$|S| = \alpha p^n$.  Then, 
$$
\LL(S)\ \geq\ (\alpha/2) \exp(-8 \alpha^{-1}c_p \log p). 
$$
\end{theorem}

\noindent {\bf Proof of the Theorem.}  
From Meshulam's theorem we know that 
if $U \subseteq \F_{p^m}$ satisfies $\EE(U) \geq \alpha/2$, and 
$m\ =\ \lceil 2c_p /\alpha \rceil$,
then $U$ contains a three-term arithmetic progression.  

Let ${\cal V}$ denote the sets of all additive subgroups of $\F_{p^n}$ of size 
$p^m$.  For our proof we will need to establish some facts about ${\cal V}$:   
First, observe that any sequence of $m$ linearly independent
vectors in $\F_{p^n}$ determines a subgroup in ${\cal V}$; however, each subgroup
has many corresponding sequences of $m$ vectors, though each subgroup has the
same number of sequences.  Now, it is easy to see that the number of sequences of 
$m$ linearly independent vectors in $\F_{p^n}$ is 
$$
(p^n - 1)(p^n - p) \cdots (p^n - p^{m-1})\ =\ \epsilon_1 p^{mn},\ {\rm where}\ 
1/2 < \epsilon_1 < 1;
$$
and, given a subgroup in ${\cal V}$ (which can also be thought of as an $\F_p$ 
vector subspace of dimension $m$), there are
$$
(p^m - 1)(p^m - p)\cdots (p^m - p^{m-1})\ =\ \epsilon_2 p^{m^2},\ {\rm where\ } 
1/2 < \epsilon_2 \leq \epsilon_1 < 1,
$$
sequences of $m$ linearly independent vectors in $\F_{p^n}$ that span this
subgroup.  So, 
$$
|{\cal V}|\ =\ \epsilon_3 p^{m(n-m)},\ {\rm where\ } 1 \leq \epsilon_3 < 2.
$$
Next, suppose that $a \in \F_{p^n}$.  We will need to know how many 
subgroups in ${\cal V}$ contain $a$:  Any such subgroup (subspace) 
can be written as ${\rm span}(a) + Z$, where ${\rm dim}(Z) = m-1$, 
and $Z \subseteq {\rm span}(a)^\perp$.
Thus, $Z$ is any $m-1$ dimensional subspace of an $n-1$ 
dimensional space; and so, from our bounds on $|{\cal V}|$, we deduce 
that there are $\epsilon_4 p^{(m-1)(n-m)}$,
$1/2 < \epsilon_4 < 1$, possibilities for $Z$, which implies that there are 
$$
\epsilon_4 p^{(m-1)(n-m)}\ =\ \epsilon_5 |{\cal V}| p^{m-n},\ {\rm where\ } 
1/2\ <\ \epsilon_5\ \leq\ 1,
$$ 
subspaces of $\F_{p^n}$ of dimension $m$ that contain $a$.  

Now, given an arithmetic progression $a,a+d,a+2d$, we note that the progression lies
in a coset $b + A$ of an additive subgroup $A$ if and only if $a \in b + A$ and $d \in A$.
Thus, if we define $T_3'(X)$ to be the number of non-trivial 
three-term arithmetic progressions belonging to a set $X$,
then the sum of the number of non-trivial arithmetic progressions lying in
$(b + A) \cap S$, over all $A \in {\cal V}$, and $b \in A^\perp$ equals
\begin{eqnarray} \label{first_double_sum}
\SUM_{A \in {\cal V}} \SUM_{b \in A^\perp} T_3'((b+A) \cap S)\ &=&\ 
\SUM_{a,a+d,a+2d \in S\atop d \neq 0} 
\SUM_{A \in {\cal V} \atop d \in A} \SUM_{b \in A^\perp \atop 
a \in b + A} 1\nonumber \\
&=&\ \SUM_{a,a+d,a+2d \in S\atop d \neq 0} 
\SUM_{A \in {\cal V} \atop d \in A} 1\nonumber \\
&\leq&\ |{\cal V}| p^{m-n} T_3'(S). 
\end{eqnarray}
We now give a lower bound on this first double sum over $A$ and $b$:  We begin with
\begin{equation} \label{AVb}
\SUM_{A \in {\cal V}} \SUM_{b \in A^\perp} |(b+A) \cap S|\ =\ |{\cal V}| |S|,
\end{equation}
which can be seen by noting that each $s \in S$ lies in exactly one coset $b+A$
of each subgroup $A \in {\cal V}$.  Now consider all the cosets 
$b + A$, $A \in {\cal V}$, such that 
\begin{equation} \label{many_cosets}
|(b + A) \cap S|\ \geq\ \alpha |A|/2. 
\end{equation}
We claim that there are more than $|{\cal V}|p^{n-m} \alpha/2$ such cosets.  To see this,
suppose there are fewer than this many cosets.  
Then, the left-most quantity in (\ref{AVb}) is at most
$$
(|{\cal V}| p^{n-m}\alpha/2) p^m + 
(|{\cal V}|p^{n-m})(\alpha |A|/2)
\ <\ |{\cal V}|\alpha p^n
\ =\ |{\cal V}| |S|,
$$ 
which would contradict (\ref{AVb}).

Thus, there are indeed more than $|{\cal V}|p^{n-m}\alpha/2$ cosets
satisfying (\ref{many_cosets}).  For each such coset $b+A$, since 
$$
|A|\ =\ p^m\ =\ p^{\lceil 2 c_p / \alpha \rceil},
$$
we deduce that $T_3'((b+A)\cap S) \geq 1$; and so, 
$$
\SUM_{A \in {\cal V}} \SUM_{b \in A^\perp} T_3'((b+A) \cap S)
\ \geq\ |{\cal V}| p^{n-m} \alpha/2. 
$$ 
Combining this with (\ref{first_double_sum}) we deduce that
$$
T_3'(S)\ \geq\ p^{2n-2m} \alpha/2\ \geq\ p^{2n} (\alpha/2) 
\exp(-8\alpha^{-1}c_p \log p).
$$
This clearly implies the theorem.

\subsection{Resumption of the Proof}

From Theorem \ref{meshulam2} and (\ref{V'_lower}) we deduce that
$$
T_3(V')\ \geq\  (\epsilon/4) \exp(-16 \epsilon^{-1} c_p \log p) |V|^2.
$$
Combining this with (\ref{T3f2}), we deduce that
$$
T_3(g)\ \leq\ T_3(f_W) - 2 \Delta p^{2n}.
$$
This, along with (\ref{fWf}) implies
$$
\LL(g)\ \leq\ \LL(f_W) - 2\Delta\ \leq\ \LL(f) - \Delta,
$$
which proves the theorem.

\end{document}